\documentclass[11pt, reqno]{amsart}

\usepackage{enumerate}
\usepackage{amsmath}%
\usepackage{amsfonts}%
\usepackage{amssymb}%
\usepackage{graphicx}
\usepackage{mathrsfs}
\usepackage{hyperref}
\usepackage{fullpage}
\usepackage{lineno}
\usepackage{setspace}
\usepackage{enumitem}
\newtheorem{theorem}{Theorem}
\theoremstyle{plain}

\newtheorem{claim}[theorem]{Claim}

\newtheorem{conjecture}[theorem]{Conjecture}

\newtheorem{definition}[theorem]{Definition}
\newtheorem{example}{Example}

\newtheorem{lemma}[theorem]{Lemma}

\newtheorem{problem}[theorem]{Problem}
\newtheorem{proposition}[theorem]{Proposition}

\numberwithin{equation}{section}
\numberwithin{theorem}{section}
\numberwithin{case}{section}

\numberwithin{subcase}{case}

\def\F{\mathcal{F}}
\def\G{\mathcal{G}}
\def\K{\mathcal{K}}
\def\h{H}
\def\G{G}

\def \a{\alpha}
\def \e{\varepsilon}
\def \r{\gamma}
\def \cP{\mathcal{P}}

\def \bfi{\mathbf{i}}

\def \bfv{\mathbf{v}}

\def\J{\mathcal{J}}

\begin{document}
\title{Covering $3$-uniform hypergraphs by vertex-disjoint tight paths}
\author{Jie Han}
\thanks{Research partially supported by Simons Foundation \#630884.}

\address{School of Mathematics and Statistics, Beijing Institute of Technology, Beijing, China, 100081}
\email{hanjie@bit.edu.cn}
\maketitle

\onehalfspace

\begin{abstract}
Let $H$ be an $n$-vertex 3-uniform hypergraph such that every pair of vertices is in at least $n/3+o(n)$ edges.
We show that $H$ contains two vertex-disjoint tight paths whose union covers the vertex set of $H$.
The quantity two here is best possible and the degree condition is asymptotically best possible.
This result also has an interpretation as the \emph{deficiency problems}, recently introduced by Nenadov, Sudakov and Wagner: every such $H$ can be made Hamiltonian by adding at most two vertices and all triples intersecting them.
\end{abstract}

\section{Introduction}

The study of Hamilton cycles is a central topic in graph theory with a long history. 
In recent years, researchers have worked on extending the classical theorem of Dirac on Hamilton cycles to hypergraphs 
and we refer to~\cite{BHS, GPW, HZ1, HZ2, RR, BMSSS1, BMSSS2, RRRSS, HZ_forbidHC} for some recent results and to \cite{KuOs14ICM, RR, zsurvey} for excellent surveys  on this topic.

In this paper we confine ourselves to 3-uniform hypergraphs (3-graphs), where each (hyper)edge contains exactly three vertices.
For a 3-graph $H$, the \emph{minimum codegree $\delta _{2} (\h)$} of $\h$ is the minimum of $\deg_{\h} (S)$ over all pairs $S$ of vertices in $\h$, where $\deg_H(S)$ is defined to be the number of edges containing $S$.
A $3$-graph $C$ is called a \emph{tight cycle} if its vertices can be ordered cyclically such that every 3 consecutive vertices in this ordering define an edge of $C$, which implies that every two consecutive edges intersect in two vertices.
We say that a $3$-graph contains a \emph{tight Hamilton cycle} if it contains a tight cycle as a spanning subgraph. 
A tight path $P$ has a sequential order of vertices $v_1v_2\dots v_{p-1}v_p$ such that every consecutive triple of vertices form an edge, where the ends of $P$ are ordered pairs $(v_2, v_1)$ and $(v_{p-1}, v_p)$.

Confirming a conjecture of Katona and Kierstead \cite{KK}, R\"odl, Ruci\'nski and Szemer\'edi \cite{RRS06, RRS11} 
determined the minimum codegree threshold for tight Hamilton cycles in 3-graphs for sufficiently large $n$, which is $\lfloor n/2\rfloor$. 
They also showed that a minimum codegree $n/2-1$ guarantees a tight Hamilton path (a spanning path).
The tightness of the results can be seen from the following example given in~\cite{KK}.
Let $V=X\dot\cup Y$, where $|X|=\lfloor n/2\rfloor$ and $|Y|=\lceil n/2\rceil$.
Let $H_0$ be a 3-graph on $V$ obtained from the complete 3-graph on $V$ by removing all triples that contain one vertex from $X$ and two vertices from $Y$.
It is straightforward to check that $\delta_2(H_0)=\min\{|X|-1, |Y|-2\}=\lceil n/2\rceil-2$.
Moreover, by construction, no tight path can connect a pair of vertices in $X$ and a pair of vertices in $Y$.
From this it is not hard to see that $H_0$ has no tight Hamilton path and adding a vertex and all triples containing it results a 3-graph with no tight Hamilton cycle.
We refer to~\cite{KK} for details.

\subsection{Main result}
A big obstruction for obtaining the tight Hamiltonicity is the `connection': even when the minimum codegree is close to $n/2$, there might be pairs of vertices that cannot be connected by a tight path.
On the other hand, the example $H_0$ above contains two vertex-disjoint tight paths whose union covers all vertices.
The aim of this paper is to show that a much weaker minimum codegree condition assures this.


\begin{theorem}\label{thm:main}
Given $\alpha>0$, there exists $n_0$ such that the following holds for $n\ge n_0$.
Let $H$ be an $n$-vertex $3$-graph with $\delta_{2}(H)\ge n/3+\alpha n$.
Then there exist two vertex-disjoint tight paths whose union covers $V(H)$.
\end{theorem}

The quantity `two' is best possible as seen by $H_0$.
The minimum codegree condition is asymptotically best possible
by the following example.
Let $V=V_0\dot\cup V_1\dot\cup V_2$, where $|V_0|=|V_1|=|V_2|=n/3$.
Let $H_1$ be the 3-graph whose edges are all triples of form $V_iV_iV_{i+1}$, $i=0,1,2$ where $V_3=V_0$.
Note that $H_1$ contains three `classes' of edges and no pair of edges from two of them can be arranged in a tight path.
Then it is easy to see that one needs three tight paths to cover $V(H_1)$.
It is also worth mentioning that $n/3$ is the asymptotical threshold for covering by \emph{a bounded number} of vertex-disjoint tight paths.
Indeed, for any $\alpha \in (0,1/3)$, the 3-graph $H_2(\alpha)$ on $V(H_2(\alpha))=X\dot\cup Y$ with $|X|=\alpha n$ and $E(H_2(\alpha))=\{e:e\cap X\neq\emptyset\}$ shows that one may need $\Omega(n)$ vertex-disjoint tight paths to cover a 3-graph $H$ with $\delta_2(H)\ge (1/3-o(1))n$.

We conjecture that $n/3$ is indeed the threshold for this property.

\begin{conjecture}\label{conj}
The minimum codegree assumption in Theorem~\ref{thm:main} can be weakened to $n/3$.
\end{conjecture}

Our problem also relates to a new type of problems, namely, `the deficiency problem', introduced very recently by Nenadov, Sudakov and Wagner~\cite{NSW_deficiency}.
Note that for Hamiltonicity, this problem has been studied previously in e.g.~\cite{Noorvash, BaOSu}.
Tailoring it to our problem, it asks for a given 3-graph $H$, what is the smallest integer $t$ such that $H*K_t^{(3)}$ contains a tight Hamilton cycle?
Here $H*K_t^{(3)}$ is a 3-graph obtained from adding $t$ vertices and all triples touching it to $H$.
The authors of~\cite{NSW_deficiency} solved the problem completely for Hamiltonicity in (2-)graphs with a given number of edges and asked for analogous results in 3-graphs. 
Theorem~\ref{thm:main} says that if $H$ satisfies that $\delta_{2}(H)\ge n/3+\e n$ for $\e >0$ and $n$ large, then $t\le 2$.
Note that the minimum-degree version of the deficiency problem for Hamiltonicity is not interesting as seen by the unbalanced complete bipartite graphs -- one needs as many vertices as to raise the minimum-degree of the resulting graph $G*K_t$ to be half of its vertices.
Thus, Theorem~\ref{thm:main} exhibits a different behavior for $k$-graphs, $k\ge 3$.

\begin{problem}
Let $H$ be a 3-graph with $\delta_2(H)\ge \alpha n$.
Determine the smallest $t=t(n, \alpha)$ such that $H*K_t^{(3)}$ contains a tight Hamilton cycle.
\end{problem}

Theorem~\ref{thm:main} shows that $t= 2$ whenever $\alpha \in (1/3, 1/2)$.
This problem makes Conjecture~\ref{conj} more appealing as an affirmative answer to the conjecture\footnote{Or a weaker statement: minimum codegree $n/3+C$ guarantees $C'$ vertex-disjoint tight paths that cover the vertex set, where $C, C'$ are absolute constants.} will give a bound on $t=t(n, \alpha)$ up to an additive constant.
Indeed, given a 3-graph $H$ with $\delta_2(H)= \alpha n$, add a set $A$ of $\beta n$ vertices and all triples touching $A$ to $H$, where $\beta$ satisfies $\frac{\alpha+\beta}{1+\beta}=\frac13$.
Thus, the resulting 3-graph has a path cover by 2 paths, and adding two more `omni' vertices results a tight Hamilton cycle, namely, $t\le \beta n+2$.
Similarly, Theorem~\ref{thm:main} shows that $t\le \beta n+o(n)$.
On the other hand, the 3-graph $H_2(\alpha)$ shows that $t\ge \beta n$.


\subsection{Proof ideas}

Our proof employs the absorbing method, which is shown to be effective on embedding spanning structures.
For example, in~\cite{RRS06}, under the minimum codegree condition $\delta_2(H)\ge (1/2+o(1))n$ it is shown that every vertex has many $v$-absorbers, a 4-vertex tight path that allows us to insert $v$ into the path as an interior vertex.
Then towards a tight Hamilton cycle, they first build an absorbing path that contains many $v$-absorbers for every vertex $v$, which can `absorb' a small but \emph{arbitrary} set of vertices and reduces the problem into finding an almost spanning tight path.
Moreover, they prove a connecting lemma: \emph{every} two pairs of vertices can be connected by a constant length tight path. 

In contrast, with a significantly weaker codegree condition, we have to look for weaker properties.

\subsubsection{Absorption}
We define our absorbers (see Section 2) for triples of vertices, and although not all triples can be absorbed, we classify the triples that can be absorbed and show that we can always partition our final leftover vertices into those triples and finish the absorption.
To classify the triples that have many absorbers, we use the lattice-based absorbing method recently developed by the author~\cite{Han14_poly}.

\subsubsection{Connection}
A pseudo-path in a $3$-graph $H$ is a sequence $(e_1,\dots, e_t)$ of not necessarily distinct edges of $H$ such that $|e_i\cap e_{i+1}|=2$ for each $i=1,\dots, t-1$.
Then a $3$-graph $H$ is connected if every two edges are connected by a pseudo-path.
The tight components of $H$ are the connected components of $H$, which are equivalence classes of edges.

For connecting, we show that \emph{among every set of three pairs of vertices, two of them can be connected by a constant length tight path}.
In fact, this is motivated by a result of Mycroft~\cite[Proposition 2]{GHM}, who proved that any $n$-vertex $3$-graph with minimum codegree $n/3$ has at most two tight components.
Inspired by this, we use the regularity method and prove that the \emph{reduced 3-graph $R$} has at most two tight components.
Then it is straightforward to show that almost every two pairs $(v_1, v_2)$, $(v_3, v_4)$ from the pair of clusters who lie in the same `component' in $R$ can be connected by a short path.
However, this only provides a connection for \emph{certain} orientation of the pairs and is not enough to prove our connecting lemma.
To see it, consider a complete 3-partite 3-graph $H$ on $V_1\cup V_2\cup V_3$, namely, every triple that meets all three clusters is an edge.
Then taking $v_1, u_1\in V_1$ and $v_2, u_2\in V_2$, there is no tight path $P$ in $H$ that connects the pairs as $v_1v_2 P u_2 u_1$.
To overcome this, we note that when the minimum codegree of $R$ is above $|R|/3$ \footnote{In fact, we only have a weaker condition: all but $o(|R|^2)$ pairs of vertices in $R$ have codegree $|R|/3$.}, every edge of $R$ lies in a copy of $K_4^-$, the unique 4-vertex 3-graph with 3 edges.
The copy of $K_4^-$ will help us to make the `turn'.
For the regularity method, we use a recent variant -- a regular slice lemma due to Allen, B\"ottcher, Cooley and Mycroft~\cite{ABCM_reg}.

Throughout the rest of the paper, by paths we mean tight paths in 3-graphs.
We introduce our absorbers in Section 2 and our connecting lemma (Lemma~\ref{lem:conn}) and path cover lemma (Lemma~\ref{lem:path}) in Section 3, followed by a proof of Theorem~\ref{thm:main}.
In Section 4, we introduce the hypergraph regularity method and the regular slice lemma (Theorem~\ref{thm:regslice}) and use them to prove Lemma~\ref{lem:conn} in Section 5 and Lemma~\ref{lem:path} in Section 6, respectively.

\section{Absorption}

In this section we give some preliminary results on the absorption part of our proof.
The following example illustrates our idea of absorbers.

\begin{example}
\label{def:abs}
Given a set of three vertices $S=\{v_1, v_2, v_3\}$, consider the following set of \emph{four} paths $P_1, P_2, P_3, P_4$
\begin{itemize}
\item for $i=1,2,3$, $P_i$ is a path $a_ib_iw_ic_id_i$ and such that $a_ib_iv_i c_i d_i$ is also a path,
\item $P_4$ is a path $u_1u_2u_3u_4$ and such that  $u_1u_2w_1w_2w_3u_3u_4$ is also a path.
\end{itemize}
\end{example}

That is, when we absorb $S$, $v_i$ will replace $w_i$ in $P_i$, $i=1,2,3$, and $w_1w_2w_3$ will be put inside $P_4$.
A known routine of the absorbing method for our problem is to show that \emph{every triple has many absorbers} and then known probabilistic arguments will produce a collection of absorbers that can absorb a small but arbitrary set of triples, which gives the existence of the absorbing path (in our problem, we may obtain a set of two paths).
Unfortunately, such a property may not hold: there might be triples that have no absorber at all and we have to classify the triples that have many absorbers.
A recent scheme to deal with such classifications is the \emph{lattice-based absorbing method} developed by the author.

To see that some triple may not have any absorber, consider the divisibility barrier: let $H_0'$ be a 3-graph with a vertex partition $V(H_0')=X\cup Y$ and the edges of $H_0'$ are all triples in $X$ and all triples that contain exactly one vertex in $X$.
When $|X|\approx |Y|$, we have $\delta_2(H_0')\approx n/2$.
Note that for any $S=\{v_1, v_2, v_3\}\subseteq Y$, since we can exchange $v_i$ and $w_i$ (as in Example~\ref{def:abs}), $w_1, w_2, w_3$ must also be in $Y$.
However, as $Y$ is an independent set, we cannot build the desired $P_4$ because any choice of $w_1w_2w_3\notin E(H_0')$.

Our actual absorbers are a little bit more complicated and allow more flexibility.
\begin{definition}
Given $S=\{v_1, v_2, v_3\}$, a family $\mathcal Q=\{P_1,\dots, P_t\}$ of vertex-disjoint tight paths is an \emph{$S$-absorber} if there exists a family of vertex-disjoint tight paths $\mathcal Q'=\{P_1',\dots, P_t'\}$ such that $V(\mathcal Q)\cup S=V(\mathcal Q')$ and $P_i'$ and $P_i$ have the same ends, for $i=1,\dots, t$, respectively.
\end{definition}

We give some notation for the lattice-based absorbing method.
Let $H$ be a $3$-graph on a vertex set $V$ with $|V|=n$.
Two (not necessarily distinct) vertices $u, v\in V$ are called \emph{$(\beta, i)$-reachable} in $H$ if there are at least $\beta n^{5i-1}$ $(5i -1)$-sets $T$ such that
\begin{itemize}
\item there exist vertex-disjoint tight paths $P_1,\dots, P_i$ of length 3 such that $V(P_1\cup\cdots\cup P_i)=T\cup \{u\}$, 
\item there exist vertex-disjoint tight paths $P_1',\dots, P_i'$ of length 3 such that $V(P_1'\cup\cdots\cup P_i')=T\cup \{v\}$, 
\item for each $j\in [i]$, $P_j$ and $P_j'$ have the same ends.
\end{itemize}
We say a vertex set $U$ is \emph{$(\beta, i)$-closed} in $H$ if any two vertices $u,v\in U$ are $(\beta, i)$-reachable in $H$.
For every $v\in V(H)$, let $\tilde{N}_{\beta, i}(v)$ be the set of vertices that are $(\beta, i)$-reachable to $v$.

We write $\alpha \ll \beta \ll \gamma$ to mean that 
it is possible to choose the positive constants
$\alpha, \beta, \gamma$ from right to left. More
precisely, there are increasing functions $f$ and $g$ such that, given
$\gamma$, whenever we choose some $\beta \leq f(\gamma)$ and $\alpha \leq g(\beta)$, the subsequent statement holds. 
Hierarchies of other lengths are defined similarly.
Given a 3-graph $H$ and a vertex $v\in V(H)$, $N_H(v)$ is defined as the collection of pairs of vertices $S\subseteq \binom{V(H)}2$ such that $S\cup \{v\}\in E(H)$.

\begin{proposition}\label{prop:Nv}
Suppose that $1/n \ll \eta \ll \alpha$. 
Let $H$ be a 3-graph with $\delta_{2}(H)\ge (1/3+\alpha)n$.
Then for any $v\in V(H)$, $|\tilde{N}_{\eta, 1}(v)| \ge  (1/3+\alpha/2)n$.
\end{proposition}

\begin{proof}
Take $\eta\ll \gamma \ll \alpha$.
Fix a vertex $v$.
For any other vertex $u\neq v$, if $|N_H(u)\cap N_H(v)|\ge \gamma n^{2}$, then by the supersaturation result (see~\cite{erdos}), there exist $\eta n^4$ copies of 4-vertex (graph) paths in $N_H(u)\cap N_H(v)$, which means that $u\in \tilde{N}_{\eta, 1}(v)$.
So if $u\notin \tilde{N}_{\eta, 1}(v)$, then  $|N_H(u)\cap N_H(v)|< \gamma n^{2}$.
By double counting, we have
\[
 (1/3+\alpha)n\cdot |N_H(v)|\le \sum_{S\in N_H(v)} \deg_{H}(S) < |\tilde{N}_{\eta, 1}(v)|\cdot |N_H(v)|+n\cdot \gamma {n}^{2}.
\]
Moreover, we have that $|N_H(v)|=\deg_H(v) \ge (n-1)(1/3+\alpha)n/2 \ge n^2/6$.
Thus, $|\tilde{N}_{\eta, 1}(v)|> (1/3+\alpha)n - {\gamma n^3}/{|N_H(v)|}\ge(1/3+\alpha/2)n$ as $\gamma \ll \alpha$.
\end{proof}

We need a ``partition lemma'' that classifies the vertices of $V(H)$ under the reachability relation.
Note that the reachability relation allows ``concatenation'', under the weakening of the constants -- that is, if for $u, v\in V(H)$, there exist at least $\mu n$ vertices $w$ such that $u$ is $(\beta_1, i_1)$-reachable to $w$ and $w$ is $(\beta_2, i_2)$-reachable to $v$, then $u$ is $(\mu \beta_1 \beta_2 - o(1), i_1+i_2)$-reachable to $v$ \footnote{Note that the $o(1)$ term exists because a pair of reachable sets for $u$ and $w$, and for $w$ and $v$ may overlap.}.


\begin{lemma}\label{lem:P}
Given $\delta>0$, and $0<\eta' \ll \delta$, there exists a constant $\beta>0$ such that the following holds for all sufficiently large $n$. 
Assume $H$ is an $n$-vertex $3$-graph such that $|\tilde{N}_{\eta', 1}(v)| \ge \delta n$ for any $v\in V(H)$. 
Then there is a partition $\cP$ of $V(H)$ into $V_1,\dots, V_r$ with $r\le 1/\delta$ such that for any $i\in [r]$, $|V_i|\ge (\delta - \eta') n$ and $V_i$ is $(\beta, 2^{\lfloor 1/\delta \rfloor-1})$-closed in $H$.
\end{lemma}

\noindent\textbf{Remark.} 
\emph{This lemma has been proved as \cite[Lemma 6.3]{HT}, under a different notion of reachability designed for the $F$-factor problem.
We remark that the form needed in this paper follows from the proof in~\cite{HT}.
Indeed, the same proof works as long as the ``reachability'' notion satisfies the following two properties:
\begin{itemize}
\item Concatenation: Let $i_1, i_2\in \mathbb N$ and $\mu, \beta_1, \beta_2\in (0,1)$. If for $u, v\in V(H)$, there exist at least $\mu n$ vertices $w$ such that $u$ is $(\beta_1, i_1)$-reachable to $w$ and $w$ is $(\beta_2, i_2)$-reachable to $v$, then $u$ is $(\mu \beta_1 \beta_2 - o(1), i_1+i_2)$-reachable to $v$.
\item Inflation: Let $i_1, i_2\in \mathbb N$ and $\beta_1\in (0,1)$. If $u$ is $(\beta_1, i_1)$-reachable to $v$ and $i_1\le i_2$, then there exists $\beta_2>0$ such that $u$ is $(\beta_2, i_2)$-reachable to $v$.
\end{itemize}
It is easy to see that our definition of reachability in this paper satisfies both points (for the second one, just append $i_2-i_1$ vertex-disjoint 5-vertex paths that are disjoint from the given reachable sets).
}

For our problem, we can take $\delta=1/3+\alpha/2>1/3$ and thus Lemma~\ref{lem:P} will return either a trivial partition or a partition of two parts, so that each part has size at least $(1/3+\alpha/3)n$ and is $(\beta, 2)$-closed in $H$.
Due to the technicality of the statement we choose not to present an absorbing lemma but rather integrate it into our proof of Theorem~\ref{thm:main}.

\section{Proof of Theorem~\ref{thm:main}}


We first present our connecting lemma and the path cover lemma whose proofs are postponed to later sections.

\begin{lemma}
[Connecting Lemma]
\label{lem:conn}

Given $\alpha>0$, there exist $\zeta_0>0$ and integer $n_0$ such that the following holds for all $\zeta <\zeta_0$ and integers $n\ge n_0$.
Let $H$ be a $3$-graph with $\delta_{2}(H)\ge (1/3+\alpha) n$.
Suppose $P_1, \dots, P_q$ are $q$ vertex-disjoint tight paths of $H$ such that $|V(P_1\cup \dots\cup P_q)|\le \zeta n$. 
Moreover, for $i=1,2$, assume that $p_i$ is a (specified) end edge of $P_i$.
Then there exist two vertex-disjoint tight paths $P_1'$, $P_2'$ such that $|V(P_1'\cup P_2')|\le \sqrt\zeta n$ and they contain $P_1,\dots, P_q$ as subpaths and contain $p_1, p_2$ as two (out of the four) ends.
\end{lemma}

Lemma~\ref{lem:conn} requires that the paths to be connected occupy a small proportion of the host graph which has a minimum codegree condition.
To use it to connect long paths in $H$, a known way is to use the trick of `reservoir': we first put aside a set $A$ of vertices chosen uniformly at random, which inherits the minimum codegree condition of $H$ even after adding a small number of other vertices in $H$ to it; then after we find the long paths, we consider $H':=H[A\cup \bigcup_i p_i]$, where $p_i$ are the ends of the paths.
So we can apply Lemma~\ref{lem:conn} as long as $|\bigcup_i p_i|\le \zeta |V(H')|$ and the connection of the ends $p_i$'s also give rise to the connection of the long paths in $H$.
One may also think of the above trick as `contracting' the long paths into 4-vertex paths.

\begin{lemma}
[Path Cover Lemma]
\label{lem:path}
Given $\alpha, \eta>0$, there exists integer $n_0$ such that the following holds for all integers $n\ge n_0$.
Let $H$ be a $3$-graph with $\delta_{2}(H)\ge (1/3+\alpha) n$.
Then there exist two vertex-disjoint tight paths $P_1$, $P_2$ such that $|V(P_1\cup P_2)|\ge n-\eta n$.
\end{lemma}

Let $\cP=\{V_1, \dots, V_\ell\}$ be a partition of $V$.
The \emph{index vector} $\bfi_{\cP}(S)\in \mathbb{Z}^\ell$ of a subset $S\subseteq V$ is the vector whose coordinates are the sizes of the intersections of $S$ with each part of $\cP$.

We recall the following Chernoff's inequality (see, e.g.,~\cite{JLR}).
For $x >0$ and a binomial random variable $X=\mathrm{Bin}(n, \zeta)$, it holds that 
\begin{equation}
\mathbb P(X\ge n\zeta + x)< e^{-x^2/(2 n\zeta + x/3)} \quad \text{and} \quad \mathbb{P}(X\le n\zeta - x)< e^{-x^2/(2 n\zeta)}. \label{eq:cher2}
\end{equation}

Now we are ready to prove Theorem~\ref{thm:main}.
The proof follows the scheme of the absorbing--reservoir method and uses Lemmas~\ref{lem:conn} and~\ref{lem:path} in the obvious way.
The additional work comes from the fact that not all triples have many absorbers.
To address this we use Lemma~\ref{lem:P} to find a partition of $V(H)$ into at most two parts, and classify the triples that do have many absorbers.
Then in the last step, we show that we can always partition the leftover vertices $A'$ into triples that have many absorbers in the absorbing paths.

\begin{proof}[Proof of Theorem~\ref{thm:main}]
Apply Lemma~\ref{lem:conn} with $\alpha/4$ in place of $\alpha$ and obtain $\zeta_0>0$.
Choose new constants $1/n\ll \eta \ll \gamma\ll \beta \ll \eta' \ll \alpha, \zeta_0$.
By Proposition~\ref{prop:Nv}, Lemma~\ref{lem:P} (applied with $\eta'\ll \delta:=1/3+\alpha$) gives a partition $\cP$ of $V(H)$ with $|\cP|=1$ or $2$ such that each part of $\cP$ has at least $(1/3+\alpha/4)n$ vertices and is $(\beta, 2)$-closed in $H$.
We first show the following claim.

\begin{claim}
For any triple $S=\{v_1, v_2, v_3\}$, if $H$ has $\alpha n^3$ edges $e$ such that $\bfi_\cP(S)=\bfi_\cP(e)$, then $H$ contains $ \beta^4 n^{34}/2$ $S$-absorbers.
\end{claim}

\begin{proof}
By the supersaturation result (see~\cite{erdos}), $H$ contains $\beta n^7$ copies of $K_{2,3,2}^{(3)}$ using these $\alpha n^3$ edges given in the claim. 
Fix a copy $K$ of such $K_{2,3,2}^{(3)}$ and take any edge $e$ from it.
Note that $\bfi_\cP(S)=\bfi_\cP(e)$.
Let $e=\{w_1, w_2, w_3\}$ such that $v_i$ and $w_i$ are $(\beta, 2)$-reachable, $i=1,2,3$.
So we can take $9$-sets $T_1, T_2, T_3$ such that for $i=1,2,3$, both $H[T_i\cup \{v_i\}]$ and $H[T_i\cup \{w_i\}]$ contain two 5-vertex paths with the same ends as stated in the definition of the reachability.
Thus for each $T_i$ there are $\beta n^9$ choices and overall there are $\beta^4 n^{34}$ choices for $K\cup T_1\cup T_2\cup T_3$.
Among them, at most $3n^{33}$ of them intersect $S$ and at most $34^2n^{33}$ of them contain repeated vertices.
Thus, there are at least $ \beta^4 n^{34}/2$ $34$-sets such that $K, T_1, T_2$ and $T_3$ are disjoint.

It remains to verify that each $K\cup T_1\cup T_2\cup T_3$ gives an $S$-absorber. 
For each $i=1,2,3$, take the two paths that span $T_i\cup \{w_i\}$ and the path $u_1u_2u_3u_4$, where $\{u_1, u_2, u_3, u_4\}=V(K)\setminus \{w_1, w_2, w_3\}$ forms a copy of $K_{1,2,1}^{(3)}$.
We claim that the family of these 7 paths is an $S$-absorber.
Indeed, for $i=1,2,3$ take the two paths that span $T_i\cup \{v_i\}$ and then take the path $u_1u_2w_1w_2w_3u_3u_4$ on $K$.
This gives a family of 7 paths which span $S\cup K\cup T_1\cup T_2\cup T_3$ and have the same ends as the family of paths mentioned above.
\end{proof}

Our main proof proceeds as the following steps.

\noindent\textbf{Build absorbing paths.}
Let $\mathcal S$ be the family of triples $S$ such that $H$ has $\alpha n^3$ edges $e$ such that $\bfi_\cP(S)=\bfi_\cP(e)$.
So the above claim says that for each $S\in \mathcal S$, $H$ contains $ \beta^4 n^{34}/2$ $S$-absorbers.
We next choose a set $\F$ of absorbers uniformly at random from $H$, that is, we select a random set $\F$ by including each $34$-set in $V(H)$ independently with probability $p:=\beta^5 n^{-33}$.
Because of \eqref{eq:cher2} (for (i) and (ii) below) and Markov's inequality (for (iii)) and the union bound, there exists a family $\mathcal F'$ satisfying the following properties:
\begin{itemize}
\item[(i)] for each triple $S\in \mathcal S$, $\mathcal F'$ contains at least $(p/2) \beta^4 n^{34}/2 = \beta^9 n/4$ $S$-absorbers;
\item[(ii)] $|\mathcal F'|\le 2p \binom n{34}\le \beta^5 n/34$;
\item[(iii)] there are at most $ 2p^2 \cdot 34\binom n{34} \binom{n}{33} \le \beta^{10} n$ pairs of overlapping members of $\mathcal F'$.
\end{itemize}
By deleting one set from each overlapping pair of members of $\F'$ and the members that are not $S$-absorbers for any triple $S\in \mathcal S$, we obtain a family $\F$ of $34$-sets such that i) $|\F|\le \beta^5 n/34$, ii) each $34$-set spans a family of 7 vertex-disjoint paths that is an $S$-absorber for some $S\in \mathcal S$, iii) for each $S\in \mathcal S$, $\F$ has at least $\beta^9 n/4 - \beta^{10} n\ge \beta^{10}n$ $S$-absorbers.

Next we use Lemma~\ref{lem:conn} with $\zeta={\beta^5}$ to connect the paths in $\F$ into two vertex-disjoint paths $P_1$ and $P_2$ such that $|V(P_1\cup P_2)|\le \sqrt{\beta^5}n\le \beta n$.
This is possible as the tight paths cover $|V(\F)|\le \beta^5 n<\zeta_0 n$ vertices.

\noindent\textbf{Choose a reservoir set $A$.}
Now we choose a random vertex set $A$ by including every vertex in $V(H)\setminus V(P_1\cup P_2)$ with probability $\gamma$.
Since $|V(P_1\cup P_2)|\le \beta n$, for any $u,v\in V(H)$, $|N_H(uv)\setminus V(P_1\cup P_2)| \ge (1/3+\alpha-\beta)n$.
By \eqref{eq:cher2} and the union bound, there exists a choice of $A$ such that $(1-2\beta)\r n\le |A|\le (1+\beta)\gamma n$ and
\begin{enumerate}[label=(\alph*)]
\item for any $u,v\in V(H)$, $|N_H(uv)\cap A| \ge (1-\beta)(1/3+\alpha-\beta)\r n \ge (1/3+\alpha/2)|A|$, \label{item:a1}
\item if $\cP=\{X, Y\}$, namely, $\cP$ has two parts, then $|A\cap X|/|A|\in (1/3+\alpha/5, 2/3-\alpha/5)$. \label{item:c}
\end{enumerate}

\noindent\textbf{Cover almost all vertices.}
Let $V'=V(H)\setminus (V(P_1\cup P_2)\cup A)$ and let $H'=H[V']$.
Since $|V(P_1\cup P_2)\cup A|\le \beta n+2\gamma n$, it holds that $\delta_2(H')\ge (1/3+\alpha/2)n$.
Then Lemma~\ref{lem:path} gives two paths $P_3$ and $P_4$ that cover all but a set $U$ of at most $\eta n$ vertices of $H'$.
We will connect the four paths to two paths by the help of $A$.
Indeed, for $i\in [4]$, we contract each $P_i$ to a 4-vertex path $\tilde P_i$ \footnote{That is, let $\tilde P_i$ denote the 4-vertex path on the end vertices of $P_i$ in order and add the possibly missing two edges to $H$. The added edges will be removed upon the completion of the connection step.} and consider $H'':=H[A\cup V(\tilde P_1\cup\cdots \cup \tilde P_4)]$.
By~\ref{item:a1}, $\delta_2(H'')\ge (1/3+\alpha/3)|V(H'')|$.
So we can apply Lemma~\ref{lem:conn} with $\zeta=16/|V(H'')| \le 16/(\r n/2)<\zeta_0$ and connect the $\tilde P_i$'s into two paths, and note that this also connects the $P_i$'s into two paths.
We take one of the paths and extend it by at most two edges so that the number of unused vertices in $A\cup U$ is a multiple of $3$. 
Denote the two paths by $Q_1$ and $Q_2$ and $A':=V(H)\setminus V(Q_1\cup Q_2)$.
It remains to absorb the vertices in $A'$.
Note that $|A'|\le |A|+|U|\le |A|+\eta n\le (1+\r)|A|$ and
\begin{equation}
\label{eq:A}
|A'| \ge |A| - \sqrt{16/|V(H'')|}|V(H'')|\ge |A| - \sqrt{16n} \ge |A| - \r|A| = (1-\r)|A|.
\end{equation}
That is, $|A'|=(1\pm \r)|A|$, and similar calculations give $|A'\cap X|=(1\pm \r)|A\cap X|$.
Together with~\ref{item:c} and $\gamma\ll \a$ we obtain $|A'\cap X|/|A'|\in (1/3+\alpha/6, 2/3-\alpha/6)$.

\noindent\textbf{Absorb the leftover.}
We first assume that $|\cP|=1$, namely, every two vertices in $V(H)$ are $(\beta, 2)$-reachable.
Since clearly $H$ contains $\alpha n^3$ edges, $\mathcal S=\binom{V(H)}3$.
In this case $Q_1$ and $Q_2$ contain $\beta^{10}n$ $S$-absorbers for every triple $S$.
As $|A'|\in 3\mathbb{N}$, $|A'|\le |A|+|U|\le 2\gamma n$ and $\gamma \ll \beta$, we can partition $A'$ arbitrarily into at most $\r n$ triples and absorb these triples one by one by their absorbers in $Q_1$ and $Q_2$.
Therefore, we obtain a path cover of $H$ by two paths.

Next assume that $\cP=\{X, Y\}$.
Let $I$ be the set of indices $\bfv\in \{(3,0), (2,1),(1,2), (0,3)\}$ such that $H$ contains $\alpha n^3$ edges $e$ with $\bfi_\cP(e)=\bfv$.
Note that we can achieve the same conclusion as in the above proof if $I=\{(3,0), (2,1),(1,2), (0,3)\}$.
We now count the edges of $H$ with different index vectors.
By $\delta_{2}(H)\ge (1/3+\alpha)n$, there are at least $\frac13\binom{|X|}2(1/3+\alpha)n > 2\alpha n^3$ edges that each contain two vertices from $X$ (recall that $|X|\ge n/3$).
This implies that one of $(3,0)$ and $(2,1)$ must be in $I$.
Similar countings derive that one of $(1,2)$ and $(2,1)$ must be in $I$, and one of $(1,2)$ and $(0,3)$ must be in $I$.

By symmetry (namely, exchange $X$ and $Y$ if necessary), it suffices to consider the following two cases.

\noindent \textbf{Case 1.} $(2,1), (1,2)\in I$.

Because $|A'\cap X|/|A'|\in (1/3+\alpha/6, 2/3-\alpha/6)$ and $|A'|\in 3\mathbb{N}$, the following system 
\[
2x^*+y^*=|A'\cap X|, \quad \text{and} \quad x^*+2y^*=|A'\cap Y|
\]
has a solution $x^*,y^*\in \mathbb N$.
So we can partition $A'$ into $x^*$ triples with index vector $(2,1)$ and $y^*$ triples with index vector $(1,2)$.
These triples can be greedily absorbed by $Q_1$ and $Q_2$ and we are done.

\noindent \textbf{Case 2.} $(2,1), (0,3)\in I$.

We need some extra work for this case.
First pick two disjoint edges $e_1$ and $e_2$ in $A'$ such that $e_1$ contains at least one vertex in $X$ and $e_2$ contains at least two vertices in $Y$.
The desired edges exist because $\delta_2(H[A'])\ge (1/3+\alpha/3)|A|$ by~\eqref{eq:A}.
Denote the specified vertex in $e_1\cap X$ by $x$ and the two specified vertices in $e_2\cap Y$ by $y_1$ and $y_2$.
Now connect the four paths $e_1, e_2, Q_1, Q_2$ into two paths $Q_1'$, $Q_2'$, so that the end with specified vertices $x$, or $y_1, y_2$ are kept as the (two out of the four) ends of $Q_1', Q_2'$.
This can be done by contracting $Q_1$ and $Q_2$ to 4-vertex paths $\tilde{Q}_1$, $\tilde{Q}_2$ and applying Lemma~\ref{lem:conn} on $H''':=H[A'\cup V(\tilde{Q}_1\cup \tilde{Q}_2)]$, because $\delta_2(H''')\ge \delta_2(H[A'])\ge (1/3+\alpha/3)|A| \ge (1/3+\alpha/4)|V(H''')|$.
Denote the set of uncovered vertices in $A'$ by $A''$.

If $|A''\cap X|$ is odd, we remove $x$, $y_1$ and $y_2$ from $Q_1'$ and $Q_2'$ (this is possible as they are at the ends) and add them to $A''$ (and we do nothing if $|A''\cap X|$ is even).
Thus $|A''\cap X|$ is even and clearly we still have $|A''|\in 3\mathbb{N}$ and $|A''\cap X|/|A''|\in (1/3, 2/3)$.
Now consider the following system
\[
2x^*=|A''\cap X|, \quad \text{and} \quad x^*+3y^*=|A''\cap Y|,
\]
which has a solution $x^*,y^*\in \mathbb N$ because $|A''\cap X|$ is even and $|A''\cap X|/|A''|\in (1/3, 2/3)$.
So we can partition $A'$ into $x^*$ triples with index vector $(2,1)$ and $y^*$ triples with index vector $(0,3)$.
These triples can be greedily absorbed and we obtain a path cover of $H$ by two paths.
\end{proof}

\section{Hypergraph regularity lemma and regular slices} \label{section:regularity}


In this section we introduce the regularity lemma and related tools we need.
The main tools needed in later proofs are the regular slice lemma (Theorem~\ref{thm:regslice}) and an extension lemma (Lemma~\ref{lem:ext}). 

\subsection{Regular complexes}
Let $\cP$ be a partition of $V$ into vertex classes $V_1, \dotsc, V_s$. A subset $S \subseteq V$ is \emph{$\cP$-partite} if $|S \cap V_i| \leq 1$ for all $1 \leq i \leq s$.
A hypergraph is \emph{$\cP$-partite} if all of its edges are $\cP$-partite, and it is \emph{$s$-partite} if it is $\cP$-partite for some partition $\cP$ with $|\cP| = s$.

A hypergraph $\h$ is a \emph{complex} if whenever $e\in E(\h)$ and $e'$ is a non-empty subset of $e$ we have that $e'\in E(\h)$.
All the complexes considered in this paper have the property that all vertices are contained in an edge.
A complex $\h$ is a \emph{$3$-complex} if all the edges of $\h$ consist of at most $3$ vertices.
The edges of size $i$ are called $i$-edges of~$\h$.
Given a $3$-complex $\h$, for all $i=1,2,3$ we denote by $\h_i$ the underlying $i$-graph of~$\h$: the vertices of $\h_i$ are those of $\h$ and the edges of $\h_i$ are the $i$-edges of~$\h$.
Given $s\ge 3$, a \emph{$(3,s)$-complex} $\h$ is an $s$-partite $3$-complex.
Given $i\le j$, an \emph{$(i,j)$-graph} is a $j$-partite $i$-graph.

Let $\h$ be a $\cP$-partite $3$-complex.
For $i \leq 3$ and $X \in \binom{\cP}{i}$, we write $\h_X$ for the subgraph of $\h_i$ induced by $\bigcup X$. 
Note that $\h_X$ is an $(i,i)$-graph.
In a similar manner we write $\h_{X^{<}}$ for the hypergraph on the vertex set $\bigcup X$, whose edge set is $\bigcup_{X' \subsetneq X} \h_{X'}$.
Note that if $\h$ is a $3$-complex and $X$ is a $3$-set, then $\h_{X^<}$ is a $(2, 3)$-complex.

Given $i\ge 2$, consider an $(i,i)$-graph $\h_{i}$ and an $(i-1,i)$-graph $\h_{i-1}$ on the same vertex set, which are $i$-partite with respect to the same partition~$\cP$.
We write $\K_i(\h_{i-1})$ for the family of all $\cP$-partite $i$-sets that form a copy of the complete $(i-1)$-graph $K_i^{i - 1}$ in~$\h_{i-1}$.
We define the \emph{density of $\h_{i}$ with respect to $\h_{i-1}$} to be
\[
d(\h_{i}|\h_{i-1})=\frac{|\K_i(\h_{i-1})\cap E(\h_{i})|}{|\K_i(\h_{i-1})|} \quad \text{if} \quad |\K_i(\h_{i-1})|>0,
\]
and $d(\h_{i}|\h_{i-1})=0$ otherwise.
More generally, if ${\bf Q}=(Q_1, \dotsc, Q_r)$ is a collection of $r$ subhypergraphs of $\h_{i-1}$, we define $\K_i({\bf Q}):=\bigcup_{j=1}^r \K_i(Q_j)$ and
\[
d(\h_{i}|{\bf Q})=\frac{|\K_i({\bf Q})\cap E(\h_{i})|}{|\K_i({\bf Q})|} \quad \text{if} \quad |\K_i({\bf Q})|>0,
\]
and $d(\h_{i}|{\bf Q})=0$ otherwise.

We say that $\h_{i}$ is \emph{$(d_i,\e,r)$-regular with respect to $\h_{i-1}$} if for all $r$-tuples~${\bf Q}$ with $|\K_i({\bf Q})|>\e |\K_i(\h_{i-1})|$ we have $d(\h_{i}|{\bf Q}) = d_i \pm \e$.
Instead of $(d_i, \e, 1)$-regularity we simply refer to \emph{$(d_i, \e)$-regularity}; we also say simply that $\h_i$ is $(\e, r)$-regular with respect to $\h_{i-1}$ if there is some $d_i>0$ for which $\h_{i}$ is $(d_i, \e, r)$-regular with respect to~$\h_{i-1}$.
Given an $i$-graph $G$ such that $V(G)\supseteq V(\h_{i-1})$, we say that~$G$ is \emph{$(d_i, \e, r)$-regular with respect to $\h_{i-1}$} if the $i$-partite subgraph of $G$ induced by the vertex classes of $\h_{i-1}$ is $(d_i, \e, r)$-regular with respect to~$\h_{i-1}$ (recall that $\h_{i-1}$ is $i$-partite).

Finally, given $s\ge 2$ and a $(2,s)$-complex $\h$ with a vertex partition $\cP$, we say that \emph{$\h$ is $(d_2,\e,r)$-regular} if for every $A \in \binom{\cP}{2}$, $\h_A$ is $(d_2, \e)$-regular with respect to $(\h_{A^<})_{1}$.
Given $s\ge 3$ and a $(3,s)$-complex $\h$ with a vertex partition $\cP$, we say that \emph{$\h$ is $(d, d_2,\e_3,\e,r)$-regular} if:
\begin{enumerate}[label={\rm (\roman*)}]
\item for every $A \in \binom{\cP}{2}$, $\h_A$ is $(d_2, \e)$-regular with respect to $(\h_{A^<})_{1}$ or $d(H_A|(H_{A<})_1)=0$, and 
\item for every $A \in \binom{\cP}{3}$, $\h_A$ is $(d, \e_3, r)$-regular with respect to $(\h_{A^<})_{2}$ or $d(H_A|(H_{A<})_2)=0$.
\end{enumerate} 
Note that by the Dense Counting Lemma (see~e.g.~\cite[Theorem 6.5]{KRS02}), a $(d, d_2,\e_3,\e,r)$-regular $(3,3)$-complex with $n$ vertices in each part and at least one 3-edge has at least $(d d_2^3/2) n^3$ $3$-edges.

We need the following lemma which states that the restriction of regular complexes to a sufficiently large set of vertices is still regular.

\begin{lemma}[Restriction Lemma,~\cite{KMO}, Lemma 4.1]
\label{lem:reg_res}
Let $s, r, m$ be positive integers and $\alpha, d_2, d, \e, \e_3>0$ such that
\[
1/m \ll 1/r, \e \le \min\{\e, d_2\}\le \e_3 \ll \alpha \ll d, 1/s.
\] 
Let $H$ be a $(d, d_2, \e_3, \e, r)$-regular $(3,s)$-complex with vertex classes $V_1, \dots, V_s$ of size $m$.
For each $i$ let $V_i'\subseteq V_i$ be a set of size at least $\alpha m$.
Then the restriction $H'=H[V_1'\cup \cdots\cup V_s']$ of $H$ to $V_1'\cup \cdots\cup V_s'$ is $(d,d_2,\sqrt{\e_3}, \sqrt\e, r)$-regular.
\end{lemma}

\subsection{Statement of the regular slice lemma}

In this section we state the version of the regularity lemma (Theorem~\ref{thm:regslice}) due to Allen, B\"ottcher, Cooley and Mycroft~\cite{ABCM_reg}, which they call the \emph{regular slice lemma}.
For most of notation in this subsection we follow those from~\cite{ABCM_reg} (with some simplification because we only focus on 3-complexes).
A similar lemma was previously applied by Haxell, \L{}uczak, Peng, R\"odl, Ruci\'{n}ski and Skokan~\cite{HLPRRS}.
This lemma says that all $3$-graphs $G$ admit a regular slice $\J$, which is a regular multipartite $2$-complex whose vertex classes have equal size such that $G$ is regular with respect to~$\J$.

Let $t_0, t_1 \in \mathbb{N}$ and $\e > 0$. 
Following~\cite{ABCM_reg}, we say that a $2$-complex $\J$ is \emph{$(t_0, t_1, \e)$-equitable} if it has the following two properties: 
\begin{enumerate}[label={\rm (\roman*)}]
\item There exists a partition $\cP$ of $V(\J)$ into $t$ parts of equal size, for some $t_0 \leq t \leq t_1$, such that $\J$ is $\cP$-partite.
We refer to $\cP$ as the \emph{ground partition} of $\J$, and to the parts of $\cP$ as the \emph{clusters} of~$\J$.
\item There exists $d_2\ge 1/t_1$ with $1/d_2 \in \mathbb{N}$, and the $2$-complex $\J$ is $(d_2, \e, 1)$-regular.
\end{enumerate}
Let $X \in \binom{\cP}{3}$.
We write $\hat{\J}_X$ for the $(2, 3)$-graph~$(\J_{X^<})_{2}$.
A $3$-graph~$G$ on $V(\J)$ is \emph{$(\e_3, r)$-regular with respect to $\hat{\J}_X$} if there exists some $d$ such that~$G$ is $(d, \e_3, r)$-regular with respect to~$\hat{\J}_X$.
We also write $d^\ast_{\J, G}(X)$ for the density of $G$ with respect to~$\hat{\J}_X$, or simply $d^\ast(X)$ if $\J$ and $G$ are clear from the context.
Now we are ready to state the definition of a regular slice from~\cite{ABCM_reg} (for 3-complexes).

\begin{definition}[Regular slice]
\cite{ABCM_reg}
Given $\e, \e_3 > 0$, $r, t_0, t_1 \in \mathbb{N}$, a $3$-graph~$G$ and a $2$-complex $\J$ on $V(G)$, we call $\J$ a \emph{$(t_0, t_1, \e, \e_3, r)$-regular slice for $G$} if $\J$ is $(t_0, t_1, \e)$-equitable and $G$ is $(\e_3, r)$-regular with respect to all but at most $\e_3 \binom{t}{3}$ of the triples of clusters of $\J$, where $t$ is the number of clusters of~$\J$.
\end{definition}

Given a regular slice $\J$ for a $3$-graph $G$, we keep track of the relative densities $d^\ast(X)$ for triples $X$ of clusters of $\J$, which is done via a weighted $3$-graph.

\begin{definition}[Weighted reduced 3-graph]\cite{ABCM_reg}
Given a $3$-graph $G$ and a $(t_0, t_1, \e)$-equitable $2$-complex~$\J$ on $V(G)$, we let $R_{\J}(G)$ be the complete weighted $3$-graph whose vertices are the clusters of $\J$, and where each edge $X$ is given weight~$d^\ast(X)$. When $\J$ is clear from the context we write $R(G)$ instead of~$R_{\J}(G)$.
\end{definition}

The regular slice lemma (Theorem~\ref{thm:regslice}) guarantees the existence of a regular slice~$\J$ with respect to which $R(G)$ resembles $G$ in various senses.
In particular, $R(G)$ inherits the codegree condition of $G$ in the following sense.
Let $G$ be a $3$-graph on $n$ vertices.
Given a set $S \in \binom{V(G)}{2}$, 
the \emph{relative degree $\overline{\deg}(S; G)$ of $S$ with respect to $G$} is defined to be 
\[ 
\overline{\deg}(S; G) = \frac{\deg_G(S)}{n - 2}, 
\]
i.e., $\overline{\deg}(S; G)$ is the proportion of triples of vertices in $G$ extending~$S$ which are in fact edges of~$G$.
To extend this definition to weighted $3$-graphs $G$ with weight function $d^\ast$, we define \[ \overline{\deg}(S; G) = \frac{\sum_{e \in E(G): S \subseteq e} d^{\ast}(e)}{n - 2}. \]
Finally, for a collection $\mathcal{S}$ of pairs in~$V(G)$, the \emph{mean relative degree $\overline{\deg}(\mathcal{S}; G)$ of $\mathcal{S}$ in $G$} is defined to be the mean of $\overline{\deg}(S; G)$ over all sets $S \in \mathcal{S}$.

We also need the `rooted counting' property in $G$ inherited by the regular slice $\J$.
For that we need the following definitions from~\cite{ABCM_reg}.
Given a $3$-graph $G$ and distinct `root' vertices $v_1,\dots, v_\ell$ of $G$, and a 3-graph $H$ with a specified set of distinct `root' vertices $x_1,\dots, x_\ell$, let $n_H(G; v_1,\dots, v_\ell)$ be the number of injective maps from $V(H)$ to $V(G)$ which embed $H$ in $G$ and map $x_j$ to $v_j$ for $j\in [\ell]$.
Then define
\[
d_H(G; v_1,\dots, v_\ell) : = \frac{n_H(G; v_1,\dots, v_\ell)}{\binom{v(G)-\ell}{v(H)-\ell}\cdot (v(H)-\ell)!}.
\]
Next we define $H^{skel}$ to be the 2-complex on $V(H)-\ell$ vertices which is obtained from the complex generated by the down-closure of $H$ by deleting the vertices $x_1,\dots, x_\ell$ and deleting all edges of size $3$.
Given a $(t_0, t_1, \e)$-equitable $2$-complex $\J$ on $V(G)$, define $n_H(G; v_1,\dots, v_\ell, \J)$ to be the number of labelled rooted copies of $H$ in $G$ such that each vertex of $H^{skel}$ lies in a distinct cluster of $\J$ and the image of $H^{skel}$ is in $\J$.
We also define $n_{H^{skel}}'(\J)$ to be the number of labelled copies of $H^{skel}$ in $\J$ with each vertex of $H^{skel}$ embedded in a distinct cluster of $\J$.
Then define
\[
d_H(G; v_1,\dots, v_\ell, \J) := \frac{n_H(G; v_1,\dots, v_\ell, \J)}{n_{H^{skel}}'(\J)}.
\]

We can now state the version of the regular slice lemma that we will use.

\begin{theorem}[Regular slice lemma~{\cite[Lemma 6]{ABCM_reg}}] 
\label{thm:regslice}
For all $t_0 \in \mathbb{N}$, $\e_3 > 0$ and all functions $r: \mathbb{N} \rightarrow \mathbb{N}$ and $\e: \mathbb{N} \rightarrow (0, 1]$, there exist $t_1, n_1 \in \mathbb{N}$ such that the following holds for all $n \ge n_1$ which are divisible by~$t_1!$.
Let $G$ be a $3$-graph on $n$ vertices.
Then there exists a $(t_0, t_1, \e(t_1), \e_3, r(t_1))$-regular slice $\J$ for $G$ such that, 
\begin{enumerate}
\item (Codegree) for all pairs $Y$ of clusters of $\J$, we have $\overline{\deg}(Y; R(G)) = \overline{\deg}(\J_Y; G) \pm \e_3$.
\item (Rooted counting) for each $1\le \ell\le 1/\e_3$, each 3-graph $H$ equipped with a set of distinct root vertices $x_1, \dots, x_\ell$ such that $v(H)\le 1/\e_3$, and any distinct vertices $v_1, \dots, v_\ell\in V(G)$, we have \label{item:root}
\[
| d_{H}(G; v_1, \dots, v_\ell, \J) - d_H(G; v_1, \dots, v_\ell) |<\e_3.
\]
\end{enumerate}
\end{theorem}

\subsection{The $d$-reduced $3$-graph and the extension lemma} \label{subsection:reducedgraph}

Once we have a regular slice $\J$ for a $3$-graph $G$, we would like to work within triples of clusters with respect to which $G$ is both regular and dense. 
To keep track of those tuples, we introduce the following definition.

\begin{definition}[The $d$-reduced $3$-graph]\cite{ABCM_reg}
Let $G$ be a $3$-graph and $\J$ be a $(t_0, t_1, \e, \e_3, r)$-regular slice for~$G$.
Then for $d > 0$ we define the \emph{$d$-reduced $3$-graph $R_d(G)$ of $G$} to be the $3$-graph whose vertices are the clusters of $\J$ and whose edges are all triples of clusters $X$ of $\J$ such that $G$ is $(\e_3, r)$-regular with respect to $X$ and $d^\ast(X) \ge d$.
Note that $R_d(G)$ depends on the choice of $\J$ but this will always be clear from the context.
\end{definition}

For $0 \leq \mu, \theta \leq 1$, we say that a $k$-graph $\h$ on $n$ vertices is \emph{$(\mu, \theta)$-dense} if there exists $\mathcal{S} \subseteq \binom{V(\h)}{k - 1}$ of size at most $\theta \binom{n}{k-1}$ such that, for all $S \in \binom{V(\h)}{k - 1} \setminus \mathcal{S}$, we have $\deg_{\h}(S) \ge \mu(n - k + 1)$. 
A $k$-graph $\h$ on $n$ vertices is \emph{strongly $(\mu, \theta)$-dense} if it is $(\mu, \theta)$-dense and, for all edges $e \in E(\h)$ and all $(k-1)$-sets $X \subseteq e$, $\deg_{\h}(X) \ge \mu (n - k + 1)$.
The next lemma was proved in~\cite{HLS_cycle}, which states that for regular slices $\J$ as in Theorem~\ref{thm:regslice}, the codegree conditions are also preserved by~$R_d(G)$.
Note that its original version allows $G$ to be $(\mu, \theta)$-dense as well.

%


\begin{lemma}~\cite{HLS_cycle} \label{lem:Rdeg}
Let $1/n \ll 1/t_1 \leq 1/t_0 \ll 1$ and $\mu, d, \e, \e_3 > 0$.
Suppose that $G$ is a $3$-graph on $n$ vertices such that $\delta_2(G)\ge \mu n$.
Let $\J$ be a $(t_0, t_1, \e, \e_3, r)$-regular slice for $G$ such that for all pairs $Y$ of clusters of $\J$, we have $\overline{\deg}(Y; R(G)) = \overline{\deg}(\J_Y; G) \pm \e_3$.
Then $R_d(G)$ is $(\mu - d - \e_3 - \sqrt{\e_3}, 3\sqrt{\e_3} )$-dense.
\end{lemma}

We use the following result proved in~\cite[Lemma 8.8]{HLS_cycle}.

\begin{lemma}
\label{lem:strongdense}
Let $n\ge 6$ and $0<\mu, \theta <1$.
Any $(\mu,\theta)$-dense $3$-graph $H$ contains a spanning subgraph $H'$ that is strongly $(\mu-8\theta^{1/4},\theta+\theta^{1/4})$-dense.
\end{lemma}

Suppose that $G$ is a $(3,\ell)$-complex with vertex classes $V_1, V_2, V_3$, and $H$ is a $(3,\ell)$-complex with vertex classes $X_1, X_2, X_3$. 
We say that $G$ \emph{respects the partition of $H$} if whenever $H$ contains an $i$-edge with vertices in $X_{j_1},\dots, X_{j_i}$, then there is an $i$-edge of $G$ with vertices in $V_{j_1},\dots, V_{j_i}$.
On the other hand, a labelled copy of $\h$ in $G$ is \emph{partition-respecting} if for each $i\in [\ell]$ the vertices corresponding to those in $X_i$ lie within $V_i$.
We write $|H|_{G}$ for the number of (labeled) partition-respecting copies of $H$ in $G$.

Roughly speaking, the Extension Lemma  says that if $\G'$ is an induced subcomplex of $\G$, and $\h$ is suitably regular, then almost all copies of $\G'$ in $\h$ can be extended to a large number of copies of $\G$ in $\h$. 
We use the following version from~\cite{ABCM_reg} which allows each triple of clusters have different densities.
\begin{lemma}
[Extension Lemma,~\cite{ABCM_reg}, Lemma 25]
\label{lem:ext}
Let $\ell, r, t, t', n_0$ be positive integers, where $t<t'$, and let $\beta, d_2, d, \e, \e_3$ be positive constants such that $1/d_2, 1/d\in \mathbb N$ and
\[
1/n_0 \ll 1/r, \e \ll c\ll \min\{\e_3, d_2\} \le \e_3 \ll \beta, d, 1/\ell, 1/t'.
\]
Then the following holds for all integers $n\ge n_0$.
Suppose that $H'$ is a $(3,\ell)$-complex on $t'$ vertices with vertex classes $Y_1,\dots, Y_\ell$ and let $H$ be an induced subcomplex of $H'$ on $t$ vertices.
Suppose also that $G$ is a $(3,\ell)$-complex with vertex classes $V_1,\dots, V_\ell$, all of size $n$, which respects the partition of $H'$, such that the 2-complex formed by 2-edges and 1-edges in $G$ is $(t_0, t_1, \e)$-equitable with density $d_2$.
Suppose further that for each 3-edge $e$ of $H'$ with index $A\in \binom{[\ell]}3$, the $(3,3)$-graph $G_A$ is $(d', \e_3, r)$-regular with respect to $(G_{A<})_2$ for some $d'\ge d$.
Then all but at most $\beta |H|_G$ labelled partition-respecting copies of $H$ in $G$ can extend to $c n^{t-t'}$ labelled partition-respecting copies of $H'$ in $G$.
\end{lemma}

\section{Connection}

In this section we prove Lemma~\ref{lem:conn}.
We first use the extension lemma (Lemma~\ref{lem:ext}) to prove the following result.

\begin{lemma}
[Connector]
\label{cor:conn}
Let $\ell, r, n_0$ be positive integers, and let $\beta, d_2, d, \e, \e_3$ be positive constants such that $1/d_2, 1/d\in \mathbb N$, $\beta\le d_2/18$ and
\[
1/n_0 \ll 1/r, \e \ll c\ll \min\{\e_3, d_2\} \le \e_3 \ll \beta, d.
\]
Then the following holds for all integers $n\ge n_0$.
Suppose that $G$ is a $(3,\ell)$-complex with a vertex partition $\cP=\{V_1,\dots, V_\ell\}$, each of size $n$, such that the 2-complex formed by 2-edges and 1-edges in $G$ is $(t_0, t_1, \e)$-equitable with density $d_2$.
Let $R$ be a $3$-graph on $\cP$ such that for each triple $T\in E(R)$, $G_T$ is $(d', \e_3, r)$-regular with respect to $(G_{T<})_2$ for some $d'\ge d$.
Let $S_1, S_0\in \binom{\cP}2$, $X_1, X_0\in \binom{\cP}3$ such that $S_i\subseteq X_i\in E(R)$, $i=0,1$, and there is a pseudo-path $P$ in $R$ that connects $X_1, X_0$, and $X_0$ is in a copy of $K_4^-$ in $R$.
Then for all but at most $\beta n^{4}$ pairs of labelled 2-edges $(v_1, v_2)$ in $G_{S_1}$ and $(v_3, v_4)$ in $G_{S_0}$, there exists a tight path $P$ of length at most $15+\ell^3$ with $(v_2, v_1)$ and $(v_3, v_4)$ as ends.
\end{lemma}

As mentioned in the introduction, the assumption that $X_0$ and $X_1$ can be connected by a pseudo-path guarantees that we can connect most of the 2-edges $v_1 v_2$ and $v_3 v_4$ but only under \emph{certain} orderings.
That is where we need the existence of a copy of $K_4^-$ to overcome the issue.

\begin{proof}

Without loss of generality, suppose $S_1=\{V_{a}, V_{a'}\}$ and $S_0=\{V_{b}, V_{b'}\}$, where $a, a', b, b'\in [\ell]$ are such that $a\neq a'$ and $b\neq b'$.
There are four cases for the pair of labelled edges, namely, e.g. for $(v_1, v_2)\in V_a\times V_{a'}$ or $V_{a'}\times V_a$ and $(v_3, v_4)\in V_b\times V_{b'}$ or $V_{b'}\times V_b$.
We will only show that for all but at most $\beta n^4/4$ pairs of labelled 2-edges $(v_1, v_2)\in V_a\times V_{a'}$ and $(v_3, v_4)\in V_b\times V_{b'}$, there exists a tight path $P$ of length at most $15+\ell^3$ with $(v_2, v_1)$ and $(v_3, v_4)$ as ends, because the same proof also treats the other three cases.

Let $X_1 X_2 \cdots X_m X_0$ be the pseudo-path $P$ in $R$ of minimum length that connects $X_1$ and $X_0$, where $X_i\in E(R)$.
Note that $|V(P)|\le 3+m$.
First we define a sequence $S=Y_1Y_2\cdots Y_{p}$ of not necessarily distinct clusters of $\cP$ (So each $Y_i=V_j$ for some $j\in [\ell]$) such that $(V_a, V_{a'})=(Y_1, Y_2)$, $S_0=\{Y_{p-1}, Y_p\}$, $X_1=\{Y_1, Y_2, Y_3\}$, $X_0=\{Y_{p-2},Y_{p-1}, Y_p\}$ and every consecutive three clusters in the sequence form one of the edges $X_0, X_1,\dots, X_m$ as follows.
We start the sequence with $Y_1Y_2 Y_3$.
After we have arranged the clusters of $X_i$ in the sequence, say $Y_{q+1}Y_{q+2} Y_{q+3}$, if $X_i\cap X_{i+1}\neq \{Y_{q+2}, Y_{q+3}\}$, then we `wind around' $X_i$, that is, let $Y_{q+4}=Y_{q+1}$ and so on, until the last two vertices of the sequence are exactly the vertices in $X_i\cap X_{i+1}$, and then put down the last vertex of $X_{i+1}$.
After having arranged $X_0$ in the sequence, wind around at most two more times if necessary so that the last two clusters in $S$ are elements of $S_0$.
Note that each time before we insert an edge $X_i$, we may need to add at most two clusters (to wind around), which implies that $|S|\le 3|V(P)| \le 3(3+m)\le 9+\ell^3$.

So the problem is that $(V_{b}, V_{b'})$ may equal $(Y_p, Y_{p-1})$, rather than $(Y_{p-1}, Y_p)$.
In this case we use the copy of $K_4^-$ to make the `turn'.
Indeed, assume that the clusters for the $K_4^-$ are $Y_{p-2}, Y_{p-1}, Y_p$ and $Y_0$.
It is straightforward to check that we can extend $S$ from the end $Y_{p-1} Y_p$ as
\begin{itemize}
\item $Y_{p-1} Y_p Y_0 Y_{p-1} Y_{p-2} Y_p Y_{p-1}$ if the missing edge of $K_4^-$ is $Y_{p-2} Y_{p} Y_0$; 
\item $Y_{p-1} Y_p Y_{p-2} Y_0 Y_{p} Y_{p-1} $ if the missing edge of $K_4^-$ is $Y_{p-2} Y_{p-1} Y_0$; 
\item $Y_{p-1} Y_p Y_{p-2} Y_0 Y_{p-1} Y_{p-2} Y_p Y_{p-1} $ if the missing edge of $K_4^-$ is $Y_{p-1} Y_{p} Y_0$.
\end{itemize}
Denote the resulting sequence by $S'=(Y_1,Y_2,\dots, Y_q)$, and thus $|S'|=q\le 15+\ell^3$.

Let $H'$ be the $(3,\ell')$-complex on distinct vertices $(x_1, x_2, \dots, x_q)$ where $\ell'=|V(P)|$ and $q=|S'|$, such that
\begin{itemize}
\item $E(H')$ is generated by the down-closure from a tight path on $(x_1, x_2, \dots, x_q)$,
\item for $i,j\in [q]$, vertices $x_i, x_j$ are in the same cluster if and only if $Y_i=Y_j$.
\end{itemize} 
Since $P$ is a pseudo-path in $R$, by the definition of $R$, $G[\bigcup V(P)]$ is a $(3,\ell')$-complex which respects the partition of $H'$ and satisfies the regularity assumptions in the lemma.
So if we let $H$ be the subcomplex of $H'$ induced on $(x_1, x_2,x_{q-1}, x_q)$, it looks like that we may apply Lemma~\ref{lem:ext} to embed $H'$.
However, this does not work because $\ell'$ might be too large, namely, we may not have, say, $\e_{3}\ll 1/\ell'$.

What we actually do is to chop $H'$ into segments and apply Lemma~\ref{lem:ext} on them separately.
More precisely, let $\ell_0=\lfloor q/9\rfloor $.
We define $H^{(1)}, \dots, H^{(\ell_0)}$ such that for $i\in [\ell_0-1]$ $H^{(i)}$ is the subcomplex of $H'$ induced on the vertices $(x_{9(i-1)+1},\dots, x_{9i+2})$, and $H^{(\ell_0)}$ is the subcomplex of $H'$ induced on $(x_{9(\ell_0-1)+1}, \dots, x_{q})$.
So each of these complexes has $11$ vertices except $H^{(\ell_0)}$ which has at most $17+2=19$ vertices.
By Lemma~\ref{lem:ext}, all but at most $\beta^2 n^{4}$ choices of pairs of labelled 2-edges $e_i\in Y_{9(i-1)+1}\times Y_{9(i-1)+2}$, $e_{i+1}\in Y_{s_i-1}\times Y_{s_i}$, where $s_i=9i+2$ for $i\in [\ell_0-1]$ and $s_{\ell_0}=q$ can be connected by at least $c n^{t_i}$ tight paths in $G$, for each $i\in [\ell_0-1]$, where $t_1=\cdots =t_{\ell_0-1}=7$ and $t_{\ell_0}=q-9(\ell_0-1)-4$.

We claim that for pairs of labelled 2-edges $e\in Y_1\times Y_2$ and $e'\in Y_{p-1}\times Y_p$, if $e$ can be connected to all but at most $8\beta n^2$ edges in $Y_{10}\times Y_{11}$ by $c n^{t_1}$ paths and $e'$ can be connected to all but at most $8\beta n^2$ edges in $Y_{9(\ell_0-1)+1} \times Y_{9(\ell_0-1)+2}$ by $c n^{t_{\ell_0}}$ paths, then $e$ and $e'$ can be connected by a desired path as stated in the lemma.
This clearly finishes the proof as the number of pairs of edges violating the properties is at most $(\beta n^2/8) \cdot n^2 + (\beta n^2/8) \cdot n^2 = \beta n^{4}/4$.
Now we prove the claim. 
Indeed, for each $j\in \{2, \dots, \ell_0-1\}$, we will only consider the $2$-edges $e_j\in Y_{9(j-1)+1}\times Y_{9(j-1)+2}$ that can be connected to all but at most $8\beta n^2$ edges in $Y_{9j+1}\times Y_{9j +2}$.
Thus, we can pick labelled 2-edges $e_2, e_3,\dots, e_{\ell_0-1}, e_{\ell_0}$ greedily so that each consecutive pair of 2-edges can be connected by at least $c n^{t_i}$ paths: indeed, when choosing $e_2, e_3,\dots, e_{\ell_0-1}$ we have at least $(d_2-\e) n^2 - 8\beta n^2 - \beta n^2/8 > \beta n^2$ choices\footnote{Recall that the $(d_2,\e)$-regularity implies the existence of $(d_2-\e) n^2$ edges. Among them, at most $8\beta n^2$ are not well connected to the previous 2-edge we chose, and at most $\beta n^2/8$ of them are not well connected to the next 2-edge to be chosen.}, and we have at least $(d_2-\e) n^2 - 2\cdot 8\beta n^2> \beta n^2$ choices for $e_{\ell_0}$.
Together with the choices for the internal vertices that connect these $e_i$'s, we have at least $\beta^{\ell_0-1} c^{\ell_0} n^{p-4}$ such candidates, of which at most $p^2n^{p-5}+4 n^{p-5}< \beta^{\ell_0-1} c^{\ell_0} n^{p-4}$ can include repeated vertices or intersect $e$ or $e'$.
So we conclude the existence of the desired path.
\end{proof}

%
%
%
%

%
The following lemma strengthens a result of Mycroft slightly, and actually follows from the same proof.
We include its (short) proof for completeness.

\begin{lemma}
\label{lem:comp}
Let $\theta\in (0,1)$.
Let $H$ be an $n$-vertex $3$-graph which is strongly $(1/3, \theta)$-dense.
Then $H$ has at most two tight components.
\end{lemma}

\begin{proof}
We regard the tight components of $H$ as an edge coloring of $H$, namely, all edges in a tight component share the same color.
Consider an edge-coloring of $K_n$, where an edge $uv$ gets the color from any $3$-edge in $H$ that contains $uv$.
This coloring is well-defined as all 3-edges containing $uv$ are in the same tight component and thus have the same color (note that an edge $uv$ may receive no color, if $\deg_H(uv)=0$).

Given a vertex $v$ and a color $c$, let $N_c(v)$ be the set of vertices that are connected to $v$ by an edge of color $c$.
Note that if $uv$ is colored with color $c$ then it has degree $n/3$ in $H$, which implies that it is adjacent to another $n/3$ edges from each of $u$ and $v$ all of color $c$.
This implies that if $uv$ and $vw$ have different colors $c_1, c_2$, then we have $|N_{c_1}(v)|, |N_{c_2}(v)|\ge n/3$ and $N_{c_1}(v)\cap N_{c_2}(v)=\emptyset$.
Therefore, there is no vertex $v$ that sees three colors, and a similar argument shows that there are no three-colored triangles.

We may assume that $H$ has three tight components, with colors $r$, $b$ and $g$. 
Since each color class of $E(K_n)$ contains a star of size $n/3$, there is a vertex, say $v$, that sees two colors, say, $r$ and $b$.
Note that $|N_r(v)|, |N_b(v)| \ge n/3$.
Since edges of $K_n$ of color $g$ contain a star of size $n/3$ and $N_{r}(v)\cap N_b(v)=\emptyset$, there is an edge $uw$ of color $g$ such that $u\in N_r(v)\cup N_b(v)$.
Note that as there is no three-colored triangle, there is no edge of color $g$ between $N_{r}(v)$ and $N_{b}(v)$.
So without loss of generality, assume that $u\in N_r(v)$ and $w\in R:=V(H)\setminus (N_{b}(v)\cup \{v\})$.
Since $|N_b(v)| \ge n/3$, we have $|R|\le |V(H)|-|N_b(v)|-1<2n/3$.
Note that there exist $n/3$ vertices $x$ which form an 3-edge of color $r$ with $uv$.
Moreover, as for such $x$, $xv$ also has color $r$, we infer $x\in N_{r}(v)\subseteq R$.
These together imply that $|N_{r}(u)\cap R|\ge n/3$.
Moreover, because $N_g(u)\subseteq R$, $|N_g(u)|\ge n/3$ and $N_{r}(u)\cap N_g(u)=\emptyset$, we derive that $|R|\ge 2n/3$, a contradiction.
\end{proof}

Now we can prove our connecting lemma.
The idea is to host the ends of all paths in the regular partition and then connect the pairs that lie in the same tight component and have the `correct' (labelled) ends required by Lemma~\ref{cor:conn}.
We remark that the technical restriction on the end edges is needed for the absorption in the proof of Theorem~\ref{thm:main}.


\begin{proof}
[Proof of Lemma~\ref{lem:conn}]
Choose constants
\[
1/n_0 \ll \zeta_0 \ll 1/t_1\ll 1/r, \e \ll c\ll \min\{\e_3, d_2\} \le \e_3 \ll d\ll \theta \ll \alpha
\]
and suppose $n\ge n_0+t_1!$ and $\zeta\le \zeta_0$.
Let $\beta = d_2^2/400$.
Let $H'$ be an induced subgraph of $H$ on $n'$ vertices such that $n'\ge n-t_1!$, $t_1!\mid n'$ and $V(P_1\cup \dots\cup P_q)\subseteq V(H')$.
We then will focus on $H'$ and note that $\delta_2(H')\ge (1/3+\alpha)n-t_1!\ge (1/3+\alpha-\theta)n$ as $1/n\le 1/n_0\ll 1/t_1\ll \theta$.
Then Theorem~\ref{thm:regslice} applied with $\e(t_1)=\e$ and $r(t_1)=r$ gives a $(t_0, t_1, \e, \e_3, r)$-regular slice $\J$ for $H'$.
Let $R_{d}(H')$ be the $d$-reduced graph which is $(1/3+\alpha-2\theta, \theta)$-dense by Lemma~\ref{lem:Rdeg}.
Then let $R$ be the strongly $(1/3+\alpha/2, 2\theta^{1/4})$-dense spanning subgraph of $R_{d}(H')$ given by Lemma~\ref{lem:strongdense}.
By Lemma~\ref{lem:comp}, $R$ has at most two tight components.
Let $\cP$ be the ground partition of $\J$ with $|\cP|=t$ and let $n_*:=n'/t$.

Let $F=x_1x_2x_3x_4$ be the 4-vertex tight path with $x_1, x_2$ as root vertices.
Since $\delta_2(H')\ge (1/3+\alpha-\theta)n$, for any $v_1, v_2\in V(H')$, by Theorem~\ref{thm:regslice}~(\ref{item:root}), we have that
\[
d_{F}(H'; v_1, v_2, \J) > d_F(H'; v_1, v_2) - \e_3 \ge \frac{(1/3+\alpha/2)^2n^2}{(n-2)(n-3)} - \e_3 \ge \frac1{9}.
\]
By the regularity, for any $X\in \binom \cP2$, it holds that $\K_2(\J_X) = (1\pm \e) d_2 n_*^2$.
Note that $F^{skel}$ is a $2$-edge (together with two singletons), and thus
\[
n_{F^{skel}}'(\J)=2\sum_{X\in \binom \cP2} \K_2(\J_X) = t(t-1) \cdot (1\pm \e) d_2n_*^2,
\]
where the factor of 2 is because $n_{F^{skel}}'(\J)$ counts \emph{labelled} copies.
These imply that $n_{F}(H'; v_1, v_2, \J)=d_{F}(H'; v_1, v_2, \J)\cdot n_{F^{skel}}'(\J)\ge \frac1{9}(1-\e)t(t-1) d_2 n_*^2$.
Since $R$ is strongly $(1/3+\alpha/2, 2\theta^{1/4})$-dense, the number of labelled copies of $F$ that are
\begin{itemize}
\item rooted at $v_1, v_2$ and
\item with $x_3, x_4$ mapped to a pair $S$ of distinct clusters of $\J$ satisfying $\deg_R(S)>0$
\end{itemize}
 is at least 
$n_{F}(H'; v_1, v_2, \J) - 2\theta^{1/4}t(t-1) (1+\e)d_2 n_*^2 \ge t(t-1) d_2 n_*^2/10$.
Therefore, there exists a pair $S:=S(v_1, v_2)$ of clusters of $\J$ such that $\deg_R(S)>0$ and $\J_S$ supports at least $d_2n_*^2/10$ labelled copies of $F$ rooted at $v_1, v_2$.
Let $H''$ be the subgraph of $H'$ which consists of the edges supported on triples of clusters in $E(R)$ only.
Since $R$ has at most two tight components, we will show that as long as there are at least three paths (which are not too long) we can connect two of them by Lemma~\ref{cor:conn}. 
Note that as we iteratively connect the paths, to guarantee the property on the end edges as stated in the lemma, it suffices to consider connecting the paths $P_1, P_2$ `at last' and when considering them with a third path we will only connect them from the end other than $p_1, p_2$.

Note also that we will use at most $q (15+t^3)$ vertices for connection and thus the collection of paths will cover at most $q (15+t^3)+\zeta n\le (16+t^3)\zeta n\le \sqrt\zeta n$ vertices, as $\zeta\le \zeta_0\ll 1/t_1\le 1/t$.
So throughout the process, by Lemma~\ref{lem:reg_res}, for each $X\in E(R)$, the restriction of $H''\cup \J$ on the set of unused vertices in $X$ is $(d', d_2, \sqrt{\e_3}, \sqrt\e, r)$-regular for some $d'\ge d$, so that we can apply Lemma~\ref{cor:conn} on the subcomplex of $H''\cup \J$ induced on the unused vertices.

Without loss of generality suppose we have paths $P_1, P_2, P_3$ and consider one end pair from each of them (but not any of $p_1, p_2$), denoted by $(v_1^i, v_2^i)$, $i=1,2,3$.
Let $S_i=S(v_1^i, v_2^i)$ be the pair of clusters defined above.
Since $\deg_R(S_i)>0$, $i=1,2,3$, take $X_i$, $i=1,2,3$, such that $S_i\subseteq X_i\in E(R)$.
As $R$ has at most two tight components, there exists $\{i,j\}\in \binom{[3]}2$ such that $X_i$ and $X_j$ are in the same tight component.
Write $X_j:=\{w_1, w_2, w_3\}$.
Since $R$ is strongly $(1/3+\alpha/2, 2\theta^{1/4})$-dense and $3(1/3+\alpha/2)n>n$, we derive that two of $N_R(w_1, w_2)$, $N_R(w_2, w_3)$ and $N_R(w_1, w_3)$ have nonempty intersection, implying the existence of a copy of $K_4^-$ containing $X_j$.
Recall that each $S_i$ hosts at least $d_2n_*^2/10$ labelled copies of $F$ rooted at $v_1^i, v_2^i$, and among them, there are at least $d_2 n_*^{2}/10 - \sqrt\zeta n\cdot n_* > d_2 n_*^{2}/20 = \sqrt{\beta} n_*^2$ such copies that do not intersect the existing paths.
As $(\sqrt{\beta} n_*^2)^2=\beta n_*^4$, there are at least $\beta n_*^4+1$ pairs of labelled copies of $F$, one rooted at $v_1^i, v_2^i$ and the other rooted at $v_1^j, v_2^j$.
If we regard the non-root vertices as labelled 2-edges, then Lemma~\ref{cor:conn} says at least one of the pairs can be connected by a tight path of length $15+t^3$, which gives the desired path connecting $(v_1^i, v_2^i)$ and $(v_1^j, v_2^j)$.
\end{proof}

\section{Path cover}

We use the following result~\cite[Lemma 4.3]{GaHa} in a slightly relaxed form\footnote{The original statement of~\cite[Lemma 4.3]{GaHa} requires that the 3-graph $H$ has no independent set of size $(2/3-o(1))n$, whose existence would imply that all pairs in the independent set has degree at most $(1/3+o(1))n$. This is indeed ruled out by our stronger degree assumption.}.

\begin{lemma} [Almost perfect matching] \label{lem:almost}
For any $\a, \theta>0$, there exist $\e_0>0$ and $n_0$ such that the following holds for $\e\le \e_0$ and $n\ge n_0$. Let $H=(V, E)$ be an $n$-vertex $3$-graph which is $(1/3+\alpha, \e)$-dense. 
Then $H$ contains a matching that covers all but at most $\theta n$ vertices of $V$.
\end{lemma}


\begin{proof}
[Proof of Lemma~\ref{lem:path}]
Apply Lemma~\ref{lem:conn} with $\alpha/3$ in place of $\alpha$ and obtain $\zeta_0>0$.
Choose constants
\[
1/n_0 \ll 1/t_1\ll 1/r, \e \ll \min\{\e_3, d_2\} \le \e_3 \ll d\ll \r \ll \theta \ll \alpha, \eta, \zeta_0
\]
and suppose $n\ge 2n_0$.
We first choose a random set $A$ of vertices by including every vertex with probability $\gamma$.
By \eqref{eq:cher2} and the union bound, there exists a choice of $A$ such that $(1-\beta)\r n\le |A|\le (1+\beta)\gamma n$ and
\begin{enumerate}[label=(\alph*)]
\item for any $u,v\in V(H)$, $|N_H(uv)\cap A| \ge (1-\beta)(1/3+\alpha)\r n \ge (1/3+\alpha/2)|A|$. \label{item:AA}
\end{enumerate}

Let $H'$ be an induced subgraph of $H-A$ on $n'$ vertices such that $n'\ge n-t_1!-|A|$ and $t_1!\mid n'$.
Note that $\delta_2(H')\ge (1/3+\alpha)n-|A|-t_1!\ge (1/3+\alpha-\theta)n$ as $1/n\le 1/n_0\ll 1/t_1\ll \theta$.
Then Theorem~\ref{thm:regslice} applied with $\e(t_1)=\e$ and $r(t_1)=r$ gives a $(t_0, t_1, \e, \e_3, r)$-regular slice $\J$ for $H'$.
Let $R_{d}(H')$ be the $d$-reduced graph which is $(1/3+\alpha/2, 3\sqrt{\e_3})$-dense by Lemma~\ref{lem:Rdeg} and the choice of the constants.
Let $\cP$ be the ground partition of $\J$ with $|\cP|=t \le t_1$ and let $n_*:=n'/t$.
By Lemma~\ref{lem:almost} applied with $\alpha/2$ in place of $\alpha$ and $3\sqrt{\e_3}$ in place of $\e$, $R_d(H')$ has a matching $M$ that covers all but $\theta t$ vertices of $R_d(H')$, and clearly $|M|\le t/3$.

Note that each edge in $M$ corresponds to a $(d', d_2, \e_3, \e, r)$-regular complex $G$ for some $d'\ge d$.
We now find a collection of vertex-disjoint paths in $G$.
Indeed, note that by Lemma~\ref{lem:reg_res}, for any subcomplex $G'$ of $G$ with $n_0\ge \theta n_*$ vertices from each cluster, $G'$ is $(d', d_2, \sqrt{\e_3}, \sqrt\e, r)$-regular.
Thus, $G'$ has at least $(dd_2^3/2)n_0^3$ $3$-edges, and by \cite[Claim 4.1]{RRS08} it contains a tight path on at least $(d d_2^3/2)n_0 \ge (dd_2^3\theta/2) n_*$ vertices.
Therefore, we can greedily construct a family of at most $6/(dd_2^3\theta)$ vertex-disjoint tight paths that together covers all but at most $3\theta n_*$ vertices.
We do the same for all edges in $M$, which altogether gives a family of at most $(t/3)(6/(dd_2^3\theta))=2t/(dd_2^3\theta)$ tight paths whose union covers all but at most $(t/3)\cdot 3\theta n_*+\theta t\cdot n_* \le 2\theta n$ vertices of $H'$.

Next we connect these paths $Q_1,\dots, Q_q$, $q\le 2t/(dd_2^3\theta)$ into two tight paths by the vertices of $A$.
To see it, for $i\in [q]$, let $p_1^i$ and  $p_2^i$  be the ends of $Q_i$ and consider $H_*:=H[A\cup \bigcup_{i\in [q]}(p_1^i\cup p_2^i)]$. 
By~\ref{item:AA}, we have that $\delta_2(H_*)\ge (1/3+\alpha/2)|A| \ge (1/3+\alpha/3)|V(H_*)|$ as $|V(H_*)|\le |A|+4q\le |A|+8t/(dd_2^3\theta)\le (1+\alpha/3)|A|$.
We regard each $p_1^ip_2^i$, $i\in [q]$ as a 4-vertex path $Q_i'$.
Because $|V(Q_1'\cup \cdots\cup Q_q')| =4q\le \zeta_0 |V(H_*)|$, we can use Lemma~\ref{lem:conn} to connect them to two paths.
This gives rise to a connection of $Q_1,\dots, Q_q$ into two tight paths $P_1, P_2$.
Note that $|V(H)\setminus V(P_1\cup P_2)|\le 2\theta n+|A|+t_1!\le \eta n$ and we are done.
\end{proof}

%

\bibliographystyle{abbrv}
\bibliography{Bibref}

\end{document}